\documentclass[11pt]{amsart}
\usepackage{amssymb}
\usepackage{amsfonts}
\usepackage{amscd}
\usepackage{amsmath}
\usepackage{mathrsfs}
\usepackage{curves}
\usepackage{epsfig}
\usepackage[pass]{geometry}
\usepackage{graphicx}
\usepackage{verbatim}
\usepackage{latexsym}
\usepackage{eucal}
\usepackage{enumerate}

\newtheorem{theorem}{Theorem}[section]
\newtheorem{lemma}[theorem]{Lemma}
\newtheorem{prop}[theorem]{Proposition}

\newtheorem{remark}[theorem]{Remark}

\def \mca {{\mathcal A}}
\def \mcb {{\mathcal B}}
\def \mcc {{\mathcal C}}
\def \mcd {{\mathcal D}}
\def \mce {{\mathcal E}}

\def \mcl {{\mathcal L}}
\def \mcm {{\mathcal M}}
\def \mcn {{\mathcal N}}

\def \mcp {{\mathcal P}}

\def \mcs {{\mathcal S}}

\def \mcv {{\mathscr V}}

\def\mcy {\mathcal{Y}}

\def \mbr {{\mathbb R}}
\def \mbs {{\mathbb S}}

\def \id {\operatorname{Id}}
\def \comp {\operatorname{comp}}
\def \loc {\operatorname{loc}}

\def \re {\operatorname{Re}}
\def \diag{\textrm{Diag}}
\def \supp {\text{supp }}

\def \beqq {\begin{equation}}
\def \eeqq {\end{equation}}
\def \WF {\text{WF}}
\def \dim {\text{dim}}
\def \bpf {\begin{proof}}
\def \epf {\end{proof}}
\def \beq {\begin{equation*}}
\def \eeq {\end{equation*}}

\def \eps {\epsilon}   
   
\def \La {\Lambda}    
   
\def \lap {\Delta}
\def \p {\partial}
\def \ha {\frac{1}{2}}

\def \tilde {\widetilde}

\def \mcx {\mathcal{X}}

\def \codim {\text{codim}}
\def \dim {\text{dim}}

\begin{document}
\title[]{Some integral geometry problems for wave equations}
\author{Yiran Wang}
\address{Yiran Wang
\newline
\indent Department of Mathematics, Emory University}
\email{yiran.wang@emory.edu}
\begin{abstract}
We consider the Cauchy problem and the source problem for normally hyperbolic operators on the Minkowski spacetime, and study the determination of solutions from their integrals along light-like geodesics. For the Cauchy problem, we give a new proof of the stable determination result obtained in Vasy and Wang \cite{VaWa}. For the source problem, we obtain stable determination for sources with space-like singularities. Our proof is based on the microlocal analysis of the normal operator of the light ray transform composed with the parametrix for strictly hyperbolic operators. 
\end{abstract}
\date{\today}
\maketitle
\section{Introduction}
Consider the $n+1$-dimensional Minkowski space  $(\mbr^{n+1}, g), n \geq 2$ where $g = -dt^2 + dx_1^2 + \cdots + dx_n^2$. Hereafter, we use $z = (z_0, z_1, \cdots, z_n) = (t, x_1, \cdots, x_n)$ for the coordinates on $\mbr^{n+1}$. Let $\square = -\p_0^2 + \sum_{j = 1}^n \p_j^2$ be the d'Alembertian where $\p_j = \frac{\p}{\p z^j}, j = 0, 1, 2, \cdots, n$. The normally hyperbolic operators on $(\mbr^{n+1}, g)$ are of the form 
\beqq\label{eq-hyper}
P(z, \p) = \square + \sum_{j = 0}^n A_j(z)\p_j + B(z)
\eeqq
where $A_j, B$ are real or complex valued smooth functions in $z$,  see e.g.\  \cite{Bar}. In this note, we study the determination of  solutions of the Cauchy problem and the source problem of \eqref{eq-hyper} from their integrals along light-like geodesics on $(\mbr^{n+1}, g)$, called the light ray transform. In addition to their own interest, these integral geometry problems arise from some inverse problems in cosmology which concern the determination of primordial gravitational perturbations from the Cosmic Microwave Background (CMB), see \cite{VaWa} for further discussions.  

To describe the light ray transform, we parametrize the future pointing light-like geodesics as follows: for $(y, \theta) \in \mcc \doteq \mbr^n\times \mbs^{n - 1}$, the light-like  geodesics from $(0, y)$ in the direction $(1, \theta)$ is given by $\gamma_{y, \theta}(s) = (s, y + s\theta),   s\in \mbr.$ 
Then the light ray transform is 
\beqq\label{eq-minlray}
L  f(y, \theta) = \int_\mbr f(s, y + s\theta) ds, \quad f\in C_0^\infty(\mbr^{n+1})
\eeqq 
It is worth mentioning that the light ray transform depends on the choice of the parametrization, see \cite[Corollary 6.2]{LOSU1}. Although light-like geodesics are preserved under conformal transformations, the light ray transform is not. 
It is also known that $L$ is injective on $C_0^\infty(\mbr^{n+1})$ see for example \cite{LOSU2} which is a result of Fourier slice theorem and the analyticity of the Fourier transform of $f$.

For $t_0 < t_1$, we denote $\mcm =  (t_0, t_1)\times \mbr^n$ and $\mcs = \{t_0\}\times \mbr^n$. For simplicity, we assume $t_0 = 0$. Consider the Cauchy problem
\beqq\label{eq-cauchy}
\begin{gathered}
P(z, \p)u(z) = 0, \text{ on } \mcm   \\
u = f_1, \quad \p_t u = f_2 \text{ on } \mcs  
\end{gathered}
\eeqq  
Our main result is the stable determination of $u$ from $L u.$
 \begin{theorem}\label{thm-main}
Let $u$ be the solution of \eqref{eq-cauchy} on $\mcm$ with Cauchy data $f_1 \in H^{s+1}(\mcs), f_2\in H^{s}(\mcs), s \geq 0$ supported in a compact set $\mcv$ of $\mcs$. Suppose that 
 the coefficients $A_j(z)$ in \eqref{eq-hyper} are real valued smooth functions. Then $f_1, f_2$ are uniquely determined by $Lu$. Furthermore, there exists $C >0$ such that 
\beqq\label{eq-stamain}
\|u\|_{H^{s+1}(\mcm)} \leq C \|(f_1, f_2)\|_{H^{s+1}(\mcs)\times H^{s}(\mcs)}\leq C \|L u\|_{H^{s +n/2 + \delta}(\mcc)}
\eeqq
where $\delta = 0$ for $n\geq 3$ and $\delta = -1/4$ for $n=2.$
\end{theorem}

The theorem for $n=3$ was proved in \cite{VaWa}. Here, the result is generalized to $n\geq 2$ and the Sobolev order in the stability estimate \eqref{eq-stamain} is improved. Despite these improvements, the main contribution of this note is to give a new proof which explores the microlocal structure of the light ray transform. It is expected that the new approach would work in more general settings. Let's recall the approach  in \cite{VaWa} and point out the differences. For the Cauchy problem \eqref{eq-cauchy}, one can use the parametrix $E$ constructed by Duistermaat and H\"ormander to represent the solution up to a smooth term. Roughly, we write $f = (f_1, f_2)$ and solution of \eqref{eq-cauchy} as $u = Ef$, thus $Lu = L  E f$.  In \cite{VaWa}, it is shown that $LE$ can be modified to an elliptic pseudo-differential operator on $\mcs$ modulo some lower order Fourier integral operators by integrating in $\theta$ variable. Then a microlocal parametrix can be constructed from which the stability estimate follows.  

In this note, we will look at the normal operator $E^\ast  L ^\ast  L E$ which seems natural to examine for an integral geometry problem. It turns out that the composition is not good as it stands. We first explain the issue in Section \ref{sec-comp1} by using a model problem. In fact, the issue is related to the microlocal structure of the normal operator $N = L^*L$. As shown in \cite{Wan1}  and reviewed in Section \ref{sec-rep},  the Schwartz kernel of $N$ is a paired Lagrangian distribution. By  judicious use of the kernel on one of the Lagrangians, we show that the composition $E^* N E$ can be slightly modified to behave well within the clean FIO calculus of Duistermaat and Guillemin, yielding a pseudo-differential operator on $\mcs$. The rest of the proof goes as in \cite{VaWa} but with less technicalities. 

We remark that in \cite{Wan1}, it is shown that for globally hyperbolic Lorentzian manifolds without conjugate points, the normal operator of the light ray transform also has a paired Lagrangian structure. The parametrix construction for the Cauchy problem for strictly hyperbolic operators works in this generality as well. Thus we believe that the new method would give stability estimates as in Theorem \ref{thm-main} and injectivity of the light ray transform for functions with sufficiently small support in such settings. These will be pursued somewhere else.   

In  Section \ref{sec-source}, we analyze the source problem from the same point of view
\beqq\label{eq-source0}
\begin{gathered}
P(z, \p) u = f, \text{ on } \mcm\\
u = 0 \text{ for } t < t_0, 
\end{gathered}
\eeqq
that is we determine $u$  on $\mcm$ from $Lu$. As pointed out however  not addressed in \cite{VaWa}, this problem arises from the inverse Sachs-Wolfe problem when the entropy perturbation cannot be ignored in Bardeen's equation.  For the source problem, there is a parametrix $E$ constructed by Melrose and Uhlmann \cite{MeUh} whose Schwartz kernel is a paired Lagrangian distribution. As for the composition $LE$, our idea is to consider $L^*LE = NE$ which turns out to be a paired Lagrangian distribution in view of a composition result of Antoniano and Uhlmann \cite{AnUh}. We will show that by considering the information on the other Lagrangian of the pair, one can stably determine $f$ when the wave front set is space-like, see Theorem \ref{thm-main1}. The case of light-like singularities is unclear. In view of the results in \cite{Wan1}, one may not be able to determine light-like singularities of $f$ in a stably way. For sources with special type of singularities such as conormal, it is possible to recover light-like singularities as in \cite{Wan1} but the result may depend on the coefficients $A_j$ in $P(z, \p)$. 

The note is organized as follows. We analyze a model problem in Section \ref{sec-mod} and   \ref{sec-comp1} where we can use oscillatory integral representations to explain the idea of the proof. Then we examine the argument from the Lagrangian distribution point of view and prove Theorem \ref{thm-main} in Section \ref{sec-rep}, \ref{sec-wavepara} and \ref{sec-comp2}. Finally, we study the source problem in  Section \ref{sec-source}.

\section{A model problem}\label{sec-mod}
We start with a model problem for which we can give an elementary proof using oscillatory integrals. Another motivation to consider a simpler example first is that through the explicit calculation, we can explain some subtlety of the problem which helps to explain the treatment for general cases.    
\begin{theorem}\label{thm-mainsim}
Let  $n\geq 3$ be an odd  integer and $s\geq 0$. 
Let $u$ be the solution of the Cauchy problem
\beqq\label{eq-scalar0}
\begin{gathered}
 \square u   = 0, \quad \text{ on }  \mcm\\ 
u = f_1, \quad \p_t u = f_2, \text{ on }  \mcs
\end{gathered}
\eeqq
where $(f_1, f_2)\in H^{s+1}(\mcs)\times H^s(\mcs)$ are supported in a compact set $\mcv\subset \mcs$. 
Then $L u$ uniquely determines $u$ on $\mcm$ and $f_1, f_2$ on $\mcs$. Moreover, we have the estimate
\beq
\|f_1\|_{H^{s+1}(\mcs)} + \|f_2\|_{H^{s}(\mcs)} \leq C\|L u\|_{H^{s+n/2}(\mcc)} 
\eeq
for $C>0$ depending on $n$ and $\mcv.$
\end{theorem}

In this section, we collect the oscillatory integral representations of the solution to the Cauchy problem and the normal operator of the light ray transform.

First, consider the solution of the Cauchy problem \eqref{eq-scalar0}. It will become convenient to consider the Cauchy problem on a larger set 
\beqq\label{eq-scalar1}
\begin{gathered}
 \square u   = 0, \quad \text{ on }    \mcn = (0, T)\times \mbr^n \\
u = f_1, \quad \p_t u = f_2, \text{ on }   \mcs
\end{gathered}
\eeqq 
where $T > t_1.$ 
Let $(\tau, \xi), \xi\in \mbr^n$ be the dual variables in $T^*\mcn$ to $(t, x), x\in \mbr^n$. 
Using Fourier transform in the $x$ variable, we get 
\beqq\label{eq-cauchysol}
\begin{gathered}
u(t, x)  
 =   (2\pi)^{-n}\int_{\mbr^n} e^{i(x\cdot \xi + t |\xi|)}  \hat h_1(\xi) d\xi +  (2\pi)^{-n}\int_{\mbr^n} e^{i(x\cdot \xi - t |\xi|)}  \hat h_2(\xi) d\xi \\
  = E_+ h_1 + E_- h_2, 
\end{gathered}
\eeqq
where 
\beq
\hat h_1 = \ha (\hat f_1 + \frac{1}{i|\xi|}\hat f_2), \ \ \hat h_2 =  \ha (\hat f_1 - \frac{1}{i|\xi|}\hat f_2).
\eeq
Here, $h_1, h_2$ are the re-parametrized Cauchy data for the Cauchy problem.  
Thus, $E_\pm$ are represented by oscillatory integrals 
\beqq\label{eq-paraE}
E_\pm f(t, x) =(2\pi)^{-n} \int_{\mbr^n}\int_{\mbr^n} e^{i((x-y)\cdot \xi \pm t |\xi|)}  f(y) dy d\xi.
\eeqq
The phase functions are  
$
\phi_\pm(t, x, y, \xi) = (x - y) \cdot \xi \pm  t |\xi|
$
and amplitude function $a(t, x, \xi) = 1.$ In particular, $E_\pm : \mce'(\mcs)\rightarrow \mcd'( \mcn)$ are Fourier integral operators with 
 canonical relations  
 \beqq\label{eq-cano0}
C_{wv}^\pm = \{(t, x, \zeta_0, \zeta'; y, \xi) \in T^*\mcn  \backslash 0\times T^*\mcs  \backslash 0: y = x -  t (\pm \xi/|\xi|), \zeta'  =  \xi, \zeta_0 = \pm |\xi|\}. 
\eeqq
Following the standard notation for Fourier integral operators see e.g.\ \cite{Dui}, we have $E_\pm \in I^{-\frac 14}(\mcn,  \mcs; C_{wv}^\pm)$.  It suffices to determine $h_1, h_2$ because we can easily find $f_1, f_2$ from
\beqq\label{eq-fh}
f_1 = h_1 + h_2, \ \ f_2 =  i  \lap^\ha (h_1 - h_2).
\eeqq

Next, consider the  light ray transform. On $\mcc = \mbr^{n}\times \mbs^{n-1}$, we use the standard product measure. Let $L^*$ be the adjoint of $L$. Consider the normal operator $N = L^* L$. It is computed in \cite[Theorem 2.1]{LOSU2} that 
\beq 
N f(t, x)  = \int_{\mbr^{n+1}} K_N(t, x, t', x') f(t', x') dt'dx'
\eeq 
where the Schwartz kernel  
\beqq\label{eq-norker}
K_N(t, x, t', x') = \frac{\delta(t - t' - |x - x'|) + \delta(t - t' + |x - x'|)}{|x - x'|^{n - 1}}
\eeqq 
In particular, $K_N$ can be written as an oscillatory integral  
\beqq \label{eq-kpseudo}
\begin{gathered}
K_N(t, x, t', x') =   \int_{\mbr^{n+1}} e^{i(t - t') \tau + i(x - x')\cdot \xi} k(\tau, \xi)  d\tau d\xi
\end{gathered}
\eeqq 
where 
\beqq\label{eq-ksym}
k(\tau, \xi) = C_n \frac{(|\xi|^2 - \tau^2)_+^{\frac{n - 3}{2}}}{|\xi|^{n-2}}, \quad C_n = 2\pi |\mbs^{n-2}|.
\eeqq
Here, for $s\in \mbr$,  $s_+^a, \re a> -1$ denotes the  distribution defined by $s_+^a = s^a$ if $s > 0$ and $s_+^a = 0$ if $s \leq 0$. Below, we denote by $\Psi^m(\mcx)$ the set of pseudo-differential operators of order $m$ on a smooth manifold $\mcx$.

 \section{The composition as oscillatory integrals}\label{sec-comp1}
Consider the determination of $h_1, h_2$ from $L u = L E_+ h_1 + L E_- h_2$. We look at the normal operator 
\beqq\label{eq-normal0}
E^* L^* L  E, \text{ where } E = E_\pm
\eeqq
There are issues about the composition as it is, and we will fine tune the operator so it becomes a pseudo-differential operator on $\mcs$. Moreover, we will show that the principal symbol is non-vanishing so the operator can be microlocally inverted.

Our idea is to modify the normal operator $N = L^*L$. Note that $u$ in \eqref{eq-scalar0} is supported on $\mcm$. We let $\chi_{[t_0, t_1]}$ be the characteristic function for $[t_0, t_1]$ on $\mbr$. Next, let $\chi\in C_0^\infty(\mbr)$ and $\supp \chi \subset (t_1, T)$. Note that $\chi_{[t_0, t_1]}\chi = 0.$ Then we consider the composition $E^*\chi L^* L \chi_{[t_0, t_1]} E$. The necessity of   the cut-off function $\chi$ is demonstrated in the next result. 
\begin{lemma}\label{lm-comp1}
For $n \geq 3$ odd, the composition $\chi L^*L \chi_{[t_0, t_1]}E_\pm \in I^{-n/2 + 1/4}(\mcn, \mcs; C^\pm_{wv})$ are elliptic Fourier integral operators.  
\end{lemma}
\bpf
We prove for $E_+$ below because the treatment for $E_-$ is identical.   The Schwartz kernel of $\chi N \chi_{[t_0, t_1]} E_+$ is 
\beq
\begin{split}
&K(t, x, z) \\
&=  (2\pi)^{-2n} \int_{\mbr^{n}} \int_{t_0}^{t_1} \int_{\mbr^{n}} \int_{\mbr^{n}} \int_{\mbr} e^{i(t - t') \tau + i(x - x')\cdot \xi} e^{i((x' - z)\cdot \eta + t'|\eta|)} \chi(t) k(\tau, \xi)  d\tau d\xi d\eta dt' dx'\\
 &=  (2\pi)^{-n}   \int_{t_0}^{t_1}   \int_{\mbr^{n}} \int_{\mbr}  e^{i(t - t') \tau + i(x - z)\cdot \xi + it'|\xi|}\chi(t) k(\tau, \xi) d\tau d\xi   dt'  
\end{split}
 \eeq
 where we integrated in $x', \eta.$ We make a change of variable $s =\tau - |\xi|$ so 
 \beq
 k(s, \xi) 
  =  C_n \frac{s_-^{\frac{n - 3}{2}}(s +2|\xi|)_+^{\frac{n-3}{2}}}{|\xi|^{n-2}}. 
 \eeq 
 Then 
  \beqq\label{eq-kker0}
 \begin{gathered}
K(t, x, z)   
  = (2\pi)^{-n}  \int_{t_0}^{t_1}   \int_{\mbr^{n}} \int_{\mbr}   e^{i(t - t') s  + i(x - z)\cdot \xi + it |\xi|}\chi(t) k(s, \xi)  d s d\xi   dt'    \\
   = (2\pi)^{-n}  \int_{t_0}^{t_1}   \int_{\mbr^{n}}    e^{i(x - z)\cdot \xi + it |\xi|} \chi(t) A(t - t', \xi) d\xi   dt'   
\end{gathered}
 \eeqq
   where $A$ is defined by 
     \beqq\label{eq-A}
 \begin{gathered}
A(\sigma, \xi)   = \int_\mbr e^{i\sigma s} k(s, \xi)ds =  \int_{-2|\xi|}^{0}   e^{i\sigma s}C_n \frac{s^{\frac{n - 3}{2}}(s +2|\xi|)^{\frac{n-3}{2}}}{|\xi|^{n-2}}d s  \\
    =C_n 2^{n-2}\int_{-1}^0 e^{2i\sigma|\xi| s}  s^{\frac{n - 3}{2}} (s  + 1)^{\frac{n-3}{2}}  ds
\end{gathered}
 \eeqq
For $n = 3$, 
     \beq
 \begin{gathered}
A(\sigma, \xi)   
   =   2 C_3 \int_{-1}^0 e^{2i\sigma|\xi| s}  ds     = C_3 \frac{1}{i\sigma } (1- e^{-2i\sigma |\xi|} )  |\xi|^{-1} 
\end{gathered}
 \eeq 
 Then 
\beqq\label{eq-tempk} 
\begin{gathered}
K(t, x, z) 
 = C_3 \int_{t_0}^{t_1}\int_{\mbr^{n}} e^{i(x - z)\cdot \xi + it |\xi|} \chi(t) \frac{1}{i (t-t') |\xi|}  d\xi  dt'\\
  + C_3\int_{t_0}^{t_1}  \int_{\mbr^{n}} e^{i(x - z)\cdot \xi - it |\xi| + 2it' |\xi|} \chi(t) \frac{-1}{i (t-t') |\xi|}   d\xi dt' 
  = K_1(t, x, z) + K_2(t, x, z)
 \end{gathered}
\eeqq 
where $K_1, K_2$ denotes the first and second integral above. For $K_1$, because $\chi(t)$ is supported away from $[t_0, t_1]$,  we see that $(t - t')^{-1}$ is integrable in $t'$. This is where we need the cut-off function! So we get 
\beq
K_1(t, x, z) =    \tilde C_3 \int_{\mbr^{n}} e^{i(x - z)\cdot \xi + it |\xi|} |\xi|^{-1} d\xi 
\eeq
with a non-vanishing constant $\tilde C_3$. This implies that $K_1(t, x, z)$ is the Schwartz kernel of an Fourier integral operator, denoted by $K_1$, associated with canonical relation $C_{wv}^+$ with a symbol of order $-1$. Note that the symbol $|\xi|^{-1}$ is singular at $\xi = 0$ but this can be removed by introducing a smooth cut-off function supported near $\xi = 0$, which amounts to changing $K_1$ by a smoothing operator.  

For $K_2$ we change variables via $\rho = t - t'$ and get 
\beq 
\begin{gathered}
K_2 = C_3\int_{t_0}^{t_1}  \int_{\mbr^{n}} e^{i(x - z)\cdot \xi - it |\xi| + 2it' |\xi|} \chi(t) \frac{-1}{i (t-t') |\xi|}   d\xi dt' \\
   = C_3\int_{t- t_0}^{t- t_1}  \int_{\mbr^{n}} e^{i(x - z)\cdot \xi + it |\xi| - 2i\rho |\xi|} \chi(t) \frac{1}{i \rho |\xi|}   d\xi d\rho 
 \end{gathered}
\eeq 
As $t - t_0> t - t_1 >0$, we can use integration by parts in $\rho$ to see that $K_2$ defines an Fourier integral operator associated with $C_{wv}^+$ with a symbol of order $-2.$ Thus $K$ in \eqref{eq-tempk} can be written as 
\beq 
\begin{gathered}
K(t, x, z)    
 =   \int_{\mbr^{n}} e^{i(x - z)\cdot \xi + it |\xi|} a(t, x, \xi)  d\xi   
 \end{gathered}
\eeq
where $a(t, x, \xi)$ is a symbol of order $-1$. The leading order term is 
\beqq\label{eq-tempa}
a_0(t, \xi) = C_3 \int_{t_0}^{t_1}  \chi(t) \frac{1}{i (t-t') |\xi|}   dt' 
\eeqq 

For $n\geq 5$ odd, we use \eqref{eq-A} and apply integration by parts  to get
\beq
\begin{gathered}
A(\sigma, \xi)  =  
 C_n \frac{2^{n-3}}{i\sigma|\xi|}  
\int_{-1}^0 e^{i\sigma 2|\xi| s}  (s^{\frac{n - 3}{2}}(s+1)^{\frac{n-3}{2}})' ds \\
 = (-1)C_n \frac{2^{n-3}}{i\sigma|\xi|}  
\int_{-1}^0 e^{i\sigma 2|\xi| s} \frac{n-3}{2} (s^{\frac{n - 3}{2}-1}(s+1)^{\frac{n-3}{2}} + s^{\frac{n - 3}{2}}(s+1)^{\frac{n-3}{2}-1}) ds
\end{gathered}
\eeq
Repeating the integration by part $\frac{n-3}{2}$ times, we get 
\beq
\begin{gathered}
A(\sigma, \xi)   = (-1)^{\frac{n-3}{2}}C_n \frac{1}{(2i \sigma |\xi|)^{\frac{n-3}{2}}} (\frac{n-3}{2})! ( \int_{-1}^0 e^{i\sigma 2|\xi|s}  (s+1)^{(n-3)/2}ds + 
  \int_{-1}^0 e^{i\sigma 2|\xi|s} s^{(n-3)/2}ds)  \\
 + \sum_{k, j \geq 1, k + j = (n-3)/2} c_{k, j} \int_{-1}^0 e^{i\sigma 2|\xi|s}  s^{(n-3)/2-k}(s+1)^{(n-3)/2-j}ds
\end{gathered}
\eeq
where $c_{k, j}$ are constants. So far, the boundary terms from integration by parts vanish. We continue with integration by parts to get 
\beq
\begin{gathered}
A(\sigma, \xi)   = (-1)^{\frac{n-3}{2}} C_n \frac{1}{(2i \sigma |\xi|)^{\frac{n-3}{2}} (2i\sigma |\xi|)} (\frac{n-3}{2})! + \sum_{k = 1}^M a_k(\sigma ) |\xi|^{-\frac{n-3}{2} -1 - k} \\
 + (-1)^{\frac{n-3}{2}} C_n \frac{1}{(2i \sigma |\xi|)^{\frac{n-3}{2}}(2i\sigma |\xi|)} (\frac{n-3}{2})!   e^{-i\sigma 2|\xi|} + \sum_{k = 1}^M b_k(\sigma ) e^{-i\sigma 2|\xi|}  |\xi|^{-\frac{n-3}{2} -1 - k}  \end{gathered}
\eeq
where $a_k, b_k$ are smooth in $\sigma$ for $\sigma \neq 0$, and $M$ is some integer depending on $n. $ We remark that for $n$ even, integration by parts will eventually lead to singular integrals. This is why we cannot deal with $n$ even at this point. We let 
\beq
\begin{gathered}
a_0(\sigma) =
b_0(\sigma) = (-1)^{\frac{n-3}{2}} C_n \frac{1}{(2i \sigma)^{\frac{n-3}{2}+1}}  (\frac{n-3}{2})!  
\end{gathered}
\eeq
 We see that 
\beq
\begin{gathered}
A(\sigma, \xi)   =  \sum_{k = 0}^M a_k(\sigma ) |\xi|^{-\frac{n-3}{2} -1 - k}  + \sum_{k = 0}^M b_k(\sigma ) e^{-i\sigma 2|\xi|}  |\xi|^{-\frac{n-3}{2} -1 - k}  
\end{gathered}
\eeq
Finally, using \eqref{eq-kker0}, we get 
\beqq\label{eq-ksym1} 
\begin{gathered}
K(t, x, z) 
 = \int_{t_0}^{t_1}  \int_{\mbr^{n}} e^{it|\xi| + i(x - z)\cdot \xi}  \sum_{k = 0}^M a_k(t - t') |\xi|^{-\frac{n-3}{2} -1 - k} d\xi dt \\
 + \int_{t_0}^{t_1} \int_{\mbr^{n}} e^{-it |\xi| + 2it'|\xi| + i(x - z)\cdot \xi} \sum_{k = 0}^M b_k(t - t')   |\xi|^{-\frac{n-3}{2} -1 - k}  d\xi dt
 \end{gathered}
\eeqq 
By the same arguments as for $n=3$, we see that $K$ is an FIO associated with $C_{wv}^+$  with a symbol of order $-(n-3)/2 - 1$.  The leading order term of the symbol is 
\beqq\label{eq-tempa1}
a_0(t, \xi) = C_n \int_{t_0}^{t_1}\chi(t)   \frac{1}{i (t-t') |\xi|^{\frac{n-3}{2} +1}}   dt' 
\eeqq 

 To summarize, for $n\geq 3$ odd, $\chi N\chi_{[t_0, t_1]} E_+$ is an FIO with canonical relation $C^+_{wv}$ of order  
\beq
-(n-3)/2 - 1 + n/2 - (2n + 1)/4 =  -n/2 + 1/4.
\eeq 
The principal symbol is clearly non-vanishing. This completes the proof. 
\epf
 
Next, we prove
\begin{lemma}\label{lm-comp2}
 For $n \geq 3$ odd, 
\begin{enumerate}
\item  $E_+^* \chi  N  \chi_{[t_0, t_1]} E_+$ and  $E_-^* \chi N  \chi_{[t_0, t_1]} E_-$  are elliptic pseudo-differential operators in  $\Psi^{-n/2 + 1/2}(\mcs)$.    
\item $E_-^*  \chi N  \chi_{[t_0, t_1]} E_+$ and  $E_+^*  \chi N   \chi_{[t_0, t_1]} E_-$ are smoothing operators on $\mcs.$
\end{enumerate}
\end{lemma}
\bpf
For (1), we consider $E_+^* \chi N  \chi_{[t_0, t_1]} E_+$. We know that the Schwartz kernel for $E_+^*$ is  
\beq
K_{E_+^*}(w, t, x) = (2\pi)^{-n} \int_{\mbr^n} e^{-i(x - w)\cdot \eta - i t|\eta|} d\eta.
\eeq
Using the notations in Lemma \ref{lm-comp1}, the kernel of $E_+^*   \chi N  \chi_{[t_0, t_1]} E_+$ is 
\beq
\begin{gathered}
K(w, z) = (2\pi)^{-2n}\int_{\mbr^n} \int_{\mbr^n} \int_{\mbr} \int_{\mbr^n} e^{-i(x - w)\cdot \eta - i t|\eta|} e^{i(x - z)\cdot \xi + it |\xi|}  a(t, \xi)  \chi(t) d\xi dtdx d\eta\\
 = (2\pi)^{-n}\int_{\mbr} \int_{\mbr^n}   e^{i(w - z)\cdot \eta}   a(t, \eta)  \chi(t) d\eta dt = (2\pi)^{-n} \int_{\mbr^n}   e^{i(w - z)\cdot \eta} c(\eta)  d\eta 
 \end{gathered}
\eeq
where $c(\eta) = \int_{\mbr} a(t, \eta)\chi(t)dt$ is a symbol of order $-\frac{n-3}{2} -1$. Thus, the composition $E_+^* \chi N  \chi_{[t_0, t_1]} E_+$ is a pseudo-differential operator of order $-n/2 + 1/2$ on $\mcs.$ 

Let's find the leading order term of $c(\eta)$, denoted by $c_0(\eta)$. For $n \geq  3$ odd, we use $a_0(t, \eta)$ in \eqref{eq-tempa} and \eqref{eq-tempa1} to get 
\beq
\begin{gathered}
c_0(\eta)  
= \int_{\mbr} \int_{t_0}^{t_1}  C_n \frac{1}{i(t - t')^{\frac{n-3}{2} + 1}}   |\eta|^{-\frac{n-3}{2} -1} \chi(t) dt' dt
\end{gathered}
\eeq
In the integral, $t > t'$ so the integrand it positive. Thus, $c_0(\eta)$ is non-zero. The proof for $E_-^* \chi N  \chi_{[t_0, t_1]} E_-$ is identical.  

For (2), let's consider $E_-^* \chi N  \chi_{[t_0, t_1]} E_+$. Then the Schwartz kernel is 
\beq
\begin{gathered}
K(w, z) = (2\pi)^{-2n}\int_{\mbr^n} \int_{\mbr^n} \int_{\mbr} \int_{\mbr^n} e^{-i(x - w)\cdot \eta + i t|\eta|} e^{i(x - z)\cdot \xi + it |\xi|}  a(t, \xi) \chi(t)   d\xi dtdx d\eta\\
 = (2\pi)^{-n}\int_{\mbr} \int_{\mbr^n}   e^{i(w - z)\cdot \eta + 2i t|\eta|}   a(t, \eta)  \chi(t)  d\eta dt = (2\pi)^{-n} \int_{\mbr^n}   e^{i(w - z)\cdot \eta} c(\eta)  d\eta 
 \end{gathered}
\eeq
Because $\chi(t)$ is supported on $t > t_1 > 0$, the integration in $t$ implies that the symbol $c(\eta)$ decays to infinite order as $|\eta|\rightarrow \infty$. So the operator is smoothing.  
\epf
 
Using the two lemmas, we can finish the proof of Theorem \ref{thm-mainsim}. 
\bpf[Proof of Theorem \ref{thm-mainsim}]
Let $u$ be the solution of \eqref{eq-scalar0}. We start with 
\beq
\chi L^*L u = \chi L^*L \chi_{[t_0, t_1]}E_+h_1 + \chi L^*L \chi_{[t_0, t_1]} E_-h_2.
\eeq
We apply $E_\pm^*$ to get
\beqq\label{eq-nop}
\begin{gathered}
E^*_+\chi N  u = E^*_+\chi N \chi_{[t_0, t_1]} E_+h_1 + E^*_+\chi N \chi_{[t_0, t_1]} E_- h_2 \\
E^*_-\chi N u  = E^*_-\chi N \chi_{[t_0, t_1]} E_+h_1 + E^*_-\chi N \chi_{[t_0, t_1]} E_- h_2
\end{gathered}
\eeqq
From Lemma \ref{lm-comp2}, $E^*_+\chi N \chi_{[t_0, t_1]} E_+, E^*_-\chi N \chi_{[t_0, t_1]} E_- \in \Psi^{-n/2 + 1/2}(\mcs)$ are elliptic pseudo-differential operators. There are parametrices $Q_\pm \in \Psi^{n/2 -1/2}(\mcs)$ such that \beq
Q_\pm E^*_\pm \chi N \chi_{[t_0, t_1]} E_\pm = \id + R_\pm
\eeq
 with $R_\pm$ smoothing operators. We also know from Lemma \ref{lm-comp2} that $E^*_+\chi N \chi_{[t_0, t_1]} E_-, $ $ E^*_-\chi N \chi_{[t_0, t_1]} E_+$ are smoothing operators. So we get from \eqref{eq-nop} that 
\beq
\begin{gathered}
Q_+E^*_+\chi N\chi_{[t_0, t_1]} u =  h_1 + R_1 h_1+ R_2 h_2 \\
Q_-E^*_-\chi N\chi_{[t_0, t_1]} u  =   h_2 + R_3 h_1 + R_4 h_2
\end{gathered}
\eeq
where $R_i, i = 1, 2, 3, 4$ are smoothing operators. 

Finally, for any $\rho\in \mbr$, we get the estimate 
\beqq\label{eq-rhoest1}
\begin{gathered}
\|h_1\|_{H^s(\mcs)} \leq C\| Q_+E^*_+\chi N\chi_{[t_0, t_1]} u\|_{H^s(\mcs)} + C_\rho \|h_1\|_{H^{s - \rho}(\mcs)} +   C_\rho \|h_2\|_{H^{s - \rho}(\mcs)}\\
\|h_2\|_{H^s(\mcs)} \leq C\| Q_-E^*_-\chi N\chi_{[t_0, t_1]} u\|_{H^s(\mcs)} + C_\rho \|h_1\|_{H^{s - \rho}(\mcs)} +   C_\rho \|h_2\|_{H^{s - \rho}(\mcs)}
\end{gathered}
\eeqq
for $C>0, C_\rho>0$. 
We know that 
\beq
\begin{gathered}
Q_\pm : H_{\comp}^s (\mcs) \rightarrow H_{\loc}^{s- n/2 + 1/2}(\mcs)
\end{gathered}
\eeq 
is bounded. For $E_\pm^*$, one can show directly using the oscillatory integral representations  or using the clean FIO calculus see Lemma \ref{lm-comp4} later,  that $E^*_\pm E_\pm \in \Psi^0(\mcs)$. Then we derive that $E^*_\pm: H_{\comp}^{s}(\mcn) \rightarrow  H_{\loc}^{s}(\mcs)$ is bounded. Finally, $L^*: H_{\comp}^{s}(\mcc)\rightarrow H_{\loc}^{s+ 1/2}(\mbr^{n+1})$ is bounded for $n\geq 3$, see Proposition \ref{prop-sobo} in Section \ref{sec-rep}. We thus conclude that    
\beq
E_\pm^* \chi L^* :   H_{\comp}^{s }(\mcc)\rightarrow H_{\loc}^{s + 1/2}(\mcn)
\eeq
is bounded. Therefore, from \eqref{eq-rhoest1}, we get  
\beqq\label{eq-rhoest2}
\begin{gathered}
\|h_1\|_{H^s(\mcs) } \leq C\| L u\|_{H^{s + n/2 - 1}(\mcc) } + C_\rho \|h_1\|_{H^{s - \rho}(\mcs) } +   C_\rho \|h_2\|_{H^{s - \rho}(\mcs) }\\
\|h_2\|_{H^s(\mcs) } \leq C\| L u\|_{H^{s + n/2-1}(\mcc) } + C_\rho \|h_1\|_{H^{s - \rho}(\mcs) } +   C_\rho \|h_2\|_{H^{s - \rho}(\mcs) }
\end{gathered}
\eeqq
Now, as shown in \cite[Theorem 8.1]{VaWa}, $L$ is injective on $L^1_{\comp}(\mbr^{n+1})$ hence on $L^2_{\comp}(\mbr^{n+1})$.  For $s\geq 0$, we know that $u$ in \eqref{eq-scalar0} belongs to $L^2_{\comp}(\mbr^{n+1})$ hence the injectivity result can be applied. One can drop the last two terms in each of the inequalities in \eqref{eq-rhoest2} by using the same argument in \cite[Theorem 1.1]{VaWa}. We get 
\beq
\begin{gathered}
\|h_1\|_{H^{s+1}(\mcs) } \leq C \| L u\|_{H^{s + n/2}(\mcc) },  \quad 
\|h_2\|_{H^{s+1}(\mcs) } \leq C\| L u\|_{H^{s + n/2}(\mcc) }  
\end{gathered}
\eeq
In terms of $f_1, f_2$ see \eqref{eq-fh}, we get
\beq
\|f_1\|_{H^{s+1}(\mcs) } + \|f_2\|_{H^{s }(\mcs) } \leq C\| L u\|_{H^{s + n/2}(\mcc) }
\eeq
This completes the proof of Theorem \ref{thm-mainsim}. 
\epf

 \section{Representation of operators}\label{sec-rep}
To understand the mechanism behind the composition in Lemma \ref{lm-comp1}, \ref{lm-comp2}, we will examine the arguments from the Lagragian distribution point of view. Note that the signature of $g$ is $(-, +, \cdots, +).$ On the dual space $\mbr^{n+1}_{(\tau, \xi)}$, we let $\Gamma^{tm}_\pm = \{(\tau, \xi)\in \mbr^{n+1}: \tau^2 > |\xi|^2, \pm \tau > 0\}$ be the set of future/past pointing time-like vectors, and $\Gamma^{tm} = \Gamma^{tm}_+\cup \Gamma^{tm}_-$. Let $\Gamma^{sp} =  \{(\tau, \xi)\in \mbr^{n+1}: \tau^2 < |\xi|^2\}$ be the set of space-like  vectors. Finally, let $\Gamma^{lt}_\pm = \{(\tau, \xi)\in \mbr^{n+1}: \tau^2 = |\xi|^2, \pm \xi_0 > 0\}$ be the set of future/past pointing light-like  vectors. We also let $\Gamma^{lt} = \Gamma_+^{lt}\cup \Gamma_-^{lt}$.   
We see that in \eqref{eq-ksym}, the symbol $k(\tau, \xi)$ is supported in $\Gamma^{sp}$,  is homogeneous of degree $-1$ in $(\tau, \xi)$ and smooth away from $\Gamma^{lt}$. Moreover,  $k(\tau, \xi) \sim \text{dist}((\tau, \xi), \Gamma^{lt})^{(n - 3)/2},$ for $(\tau, \xi)$ space-like near $\Gamma^{lt}.$ Therefore, $k(\tau, \xi)$ looks like a symbol for a pseudo-differential operator of order $-1$ with a conormal singularity at $\Gamma^{lt}$. This is an example of paired Lagrangian distribution introduced in \cite{GuUh}, as proved in \cite{Wan1} for general globally hyperbolic Lorentzian manifolds. In this section, we briefly recall the notion of paired Lagrangian distributions and   the construction for the Minkowski spacetime.  
 
%

To define paired Lagrangian distributions, we first consider the following model problem. 
Let $\tilde \mcx = \mbr^n = \mbr^k\times \mbr^{n-k}, 1\leq k\leq n-1$, and use coordinates $x = (x', x''), x'\in \mbr^k, x''\in \mbr^{n-k}$. Let $\tilde \La_0 = \{(x, \xi, x, -\xi)\in T^*(\tilde \mcx\times \tilde \mcx)\backslash 0 : \xi\neq 0\}$ be the punctured conormal bundle of $\diag$ in $T^*(\tilde \mcx \times \tilde \mcx)$, and 
\beq
\tilde \La_1 = \{(x, \xi, y, \eta)\in T^*(\tilde \mcx\times \tilde \mcx)\backslash 0: x'' = y'', \xi' = \eta' = 0, \xi'' = \eta'' \neq 0\} 
\eeq
which is the punctured conormal bundle to $\{(x, y)\in \tilde \mcx\times \tilde \mcx: x'' = y''\}$.  The two Lagrangians intersect cleanly at $\tilde \Sigma = \{(x, \xi, y, \eta)\in T^*(\tilde \mcx\times \tilde \mcx)\backslash 0 : x'' = y'', \xi'' = \eta'', x' = y', \xi' = \eta' = 0\}$ which is of codimension $k.$ 
For this model pair, the paired Lagrangian distribution $I^{p, l}(\mbr^n\times \mbr^n; \tilde \La_0, \tilde \La_1)$ consists of oscillatory integrals  (see \cite[Section 5]{DUV})
 \beqq\label{eq-upair1}
u(x, y) = \int e^{i[(x' - y')\cdot \eta' + (x'' - y'')\cdot \eta'']}b(x, y, \eta) d\eta 
\eeqq
modulo $C_0^\infty(\mbr^n\times \mbr^n)$, where $b$ satisfies the following estimates. We remark that the order here is different from that in \cite{DUV} because we work on the product space. First, in the region $|\eta'|\leq C |\eta''|, |\eta''|\geq 1$, $b$ satisfies
\beq
|(Qb)(x, y, \eta)|\leq C\langle \eta''\rangle^{p+k/2} \langle \eta'\rangle^{l - k/2}
\eeq
for all $Q$ which is a finite product of differential operators of the form $D_{\eta'}, \eta'_j D_{\eta'_m},$ $\eta''_j D_{\eta_m''}$. Second, in the region $|\eta''|\leq C |\eta'|, |\eta'|\geq 1$, $b$ satisfies the standard regularity estimate 
\beq
|(Qb)(x, y, \eta)|\leq C\langle \eta'\rangle^{p+ l}
\eeq
for all $Q$ which is a finite product of differential operators of the form $\eta_j'D_{\eta'_m}, \eta'_j D_{\eta''_m}.$ 
We use the notation $I^{p, l}(\mbr^n\times \mbr^n; \tilde \La_0, \tilde \La_1)$ to denote the space of operators $A: \mce'(\mbr^n; \Omega^\ha_{\mbr^n}) \rightarrow \mcd'(\mbr^n; \Omega^\ha_{\mbr^n})$ where $\Omega^\ha_{\mbr^n}$ denotes the line bundle of half-densities on $\mbr^n$,  whose Schwartz kernel $K_A$ is a paired Lagrangian distribution with values in $\Omega_{\mbr^n\times \mbr^n}^\ha$.  

Let $\mcx$ be a $C^\infty$ manifold of dimension $n$.  Let $\La_0, \La_1$ be two conic Lagrangian submanifold of $T^*(\mcx\times \mcx)\backslash 0$  such that $\La_0\cap \La_1$ cleanly at a codimension $k$, $1\leq k\leq 2n-1$ submanifold $\Sigma$. 
From \cite[Proposition 2.1]{GuUh}, we know that all such intersecting pairs $(\La_0, \La_1)$ are locally symplectic diffeomorphic to each other. 
Let $\chi: T^*(\mcx\times \mcx)\backslash 0\rightarrow T^*(\tilde \mcx\times \tilde \mcx)\backslash 0$ be a canonical transformation such that $\chi(\La_0) \subseteq \tilde \La_0, \chi(\La_1)\subseteq \tilde \La_1$.  Then the set of paired Lagrangian distributions $I^{p, l}(\mcx \times \mcx; \La_0, \La_1)$ are defined invariantly by conjugating elements of $I^{p, l}(\tilde \mbr^n\times \tilde \mbr^n; \tilde \La_0, \tilde \La_1)$ by Fourier integral operators with canonical relation $\chi$, see \cite{GuUh} for more details. For any $u\in I^{p, l}(\mcx\times \mcx; \La_0, \La_1)$, 
$u$ is an Fourier integral operator of order $p+l$ on $\La_0\backslash \Sigma$ and  an Fourier integral operator of order $p$ on $\La_1\backslash \Sigma$. The principal symbols satisfy certain compactbility conditions at $\Sigma$. In particular, the principal symbol of $u$ on $\La_1\backslash \Sigma$ is singular at $\Sigma.$


To see that the kernel \eqref{eq-norker} is a paired Lagrangian distribution, we can use the oscillatory integral representation \eqref{eq-kpseudo}, \eqref{eq-ksym}. 
We let $s = \tau - |\xi|$ and write  
\beqq\label{eq-ks}
\begin{gathered}
K_N(t, x, t', x') =   \int_{\mbr^{n+1}} e^{i(t - t')(s + |\xi|) + i(x - x')\cdot \xi} k(s, \xi)  ds d\xi, \\
 k(s, \xi)  =  C_n \frac{s_-^{\frac{n - 3}{2}}(s +2|\xi|)_+^{\frac{n-3}{2}}}{|\xi|^{n-2}}.
\end{gathered}
\eeqq
To see that this can be transformed to the model integral \eqref{eq-upair1}, consider the simplectic change of variables  on $T^*\mbr^{n+1}$
\beq
\tilde x = x - (t - t')\xi/|\xi|, \quad \tilde t = t - t', \quad s = s, \quad \xi = \xi.
\eeq
We can choose an Fourier integral operator with symbol of order $0$ which quantizes the symplectic change of variable to transform $K$ to
\beqq\label{eq-kker}
\begin{gathered}
K_N(\tilde t, \tilde x, t', x')   =  \int_{\mbr^{n+1}} e^{i\tilde t s + i \tilde x\cdot \xi} k(s, \xi) ds d\xi 
 \end{gathered}
\eeqq 
modulo a smooth term. 
The symbol $k(s, \xi)$ satisfies the product type estimate with $p = -n/2, l = n/2-1$. In fact, for $|\xi| \leq C|s|, |s|\geq 1$, we have 
\beq
|k(s, \xi)| \leq 
 C|s|^{-1}
\eeq
One can verify the same estimate for $Qk$ where $Q$ is the finite product of differential operators of the form 
$sD_{s}, s D_{\xi_m}$. For $|s|\leq C|\xi|, |\xi|\geq 1$, we have 
\beq
|k(s, \xi)| \leq C\frac{|s|^{\frac{n - 3}{2}} |\xi|^{\frac{n-3}{2}}}{|\xi|^{n-2}} \leq C  |\xi|^{-n/2 + 1/2} |s|^{n/2-3/2}
\eeq
 and one can verify the estimate for $Qk$ where $Q$ is the finite product of differential operators of the form $D_s, s D_{s}, \xi_j D_{\xi_m}$. So $K_N$ is a paired Lagrangian distribution. The two associated Lagrangians are  
 \beqq\label{eq-lag1}
 \La_0 = \{(t, x, \tau, \xi; t', x', \tau', \xi')\in T^*\mbr^{n+1}\backslash 0\times T^*\mbr^{n+1}\backslash 0: t = t', x = x', \tau = -\tau', \xi = -\xi'\}
 \eeqq
which is the punctured conormal bundle of the diagonal in $\mbr^{n+1}\times \mbr^{n+1}$ and 
 \beqq\label{eq-lag2}
 \begin{gathered}
 \La_1 = \{(t, x, \tau, \xi; t', x', \tau', \xi')\in T^*\mbr^{n+1}\backslash 0\times T^*\mbr^{n+1}\backslash 0: x = x' + (t - t')\xi/|\xi|\\
 \tau = \pm|\xi|,  \tau' = -\tau,  \xi' = -\xi\}.
  \end{gathered}
 \eeqq 
The two Lagrangians intersect cleanly at 
\beqq\label{eq-sigma}
\begin{gathered}
\Sigma = \{(t, x, \tau, \xi; t', x', \tau', \xi')\in T^*\mbr^{n+1}\backslash 0\times T^*\mbr^{n+1}\backslash 0: t = t', x = x', \\
\tau = -\tau', \xi = -\xi', \tau^2 = |\xi|^2\}
\end{gathered}
\eeqq
In fact, $\La_1$ is the flow out of $\Sigma$ under the Hamilton vector field $H_f$ of $f(\tau, \xi) = \ha(\tau^2 - |\xi|^2).$  

\begin{theorem}[Theorem 3.1 of \cite{Wan1}]\label{thm-mainmin0}
For the Minkowski light ray transform $L$ defined in \eqref{eq-minlray}, the Schwartz kernel of the normal operator $N = L^\ast L$ belongs to $I^{-n/2, n/2- 1}(\mbr^{n+1}\times \mbr^{n+1}; \La_0, \La_1)$, in which $\La_0, \La_1$ are two cleanly intersection Lagrangians defined in \eqref{eq-lag1}, \eqref{eq-lag2}.  The principal symbols  of $N$ on $\La_1\backslash \Sigma$ are  real valued and non-vanishing. 
\end{theorem}

The principal symbols of $K_N$ on $\La_0\backslash \La_1$ and $\La_1\backslash \La_0$ can be found explicitly, see \cite{Wan1} for details. We only need the symbol on $\La_1\backslash \La_0$ 
where the kernel $K_N \in I^{-n/2}(\mcm\times \mcm; \La_1)$. 
 To find the symbol, we can use \eqref{eq-norker} and write $K_N$ as 
\beqq\label{eq-knsym1} 
K_N(t, x, t', x') = \int_{\mbr} e^{i(t - t' - |x - x'|)\tau}  (t - t')_+^{-(n-1)}d\tau  +\int_{\mbr} e^{i(t - t' + |x - x'|)\tau} (t - t')_-^{-(n-1)}d\tau
\eeqq 
for  $t\neq t'$. This gives another oscillatory integral representation of $K_N$ with a real phase function valid for $t\neq t'.$ 
The principal symbol is non-vanishing and positive in this representation.

Using the  estimates for paired Lagrangian distributions for the flow out model 
one can derive the Sobolev estimates for $L$ and $L^*$. 
\begin{prop}[Corollary 3.2 of \cite{Wan1}]\label{prop-sobo}
The Minkowski light ray transform $L: H^s_{\comp}(\mbr^{n+1}) \rightarrow H_{\loc}^{s+ s_0/2}(\mbr^n\times \mbs^{n-1})$ and its adjoint $L^*: H_{\comp}^{s}(\mbr^n\times \mbs^{n-1}) \rightarrow H^{s+ s_0/2}_{\loc}(\mbr^{n+1}) $ are continuous 
for $s_0 = 1/2$ when $n=2$ and for $s_0 = 1$ when $n\geq 3.$
\end{prop} 
 
\section{Parametrix for the Cauchy problem}\label{sec-wavepara}
A linear differential operator $P: C^\infty(\mbr^{n+1})\rightarrow C^\infty(\mbr^{n+1})$ of second order is called {\em normally hyperbolic} if the principal symbol $\mcp(z, \zeta) \doteq \sigma(P)(z, \zeta) = g^*(\zeta, \zeta), (z, \zeta) \in T^*M$, see  \cite[page 33]{Bar}.   Note that $P(z, \p)$ in \eqref{eq-hyper} is exactly the normally hyperbolic operator on  $(\mbr^{n+1}, g)$. 
The operator is strictly hyperbolic of multiplicity one with respect to the Cauchy hypersurfaces $\mcs_t = \{t\}\times \mbr^n, t\in \mbr$, see Definition 5.1.1 of \cite{Dui}. This means that all bicharacteristic curves of $P$ are transversal to $\mcs_t$ and for $(\bar z, \bar \zeta) \in T^*\mcs_t \backslash 0$
\beq
\mcp(\bar z, \zeta) = 0, \quad \zeta|_{T_{\bar z}\mcs} = \bar \zeta
\eeq
has exactly one solution. As before,  we also use $\mcs_0 = \mcs$. It is convenient to use $D_j = -\imath \p_j, j = 0, 1, 2, \cdots, n$ in which $\imath^2 = -1.$ 
Consider the Cauchy problem 
\beqq\label{eq-cauchy1}
\begin{gathered}
P(z, D)u(z) = 0, \text{ on } \mcm \\
u = f_1, D_t u = f_2 \text{ on } \mcs.
\end{gathered}
\eeqq
We use Duistermaat-H\"ormander's parametrix construction, see e.g.\ \cite{Dui}. 
The restriction operator $\rho_0: C^\infty(\mcn) \rightarrow C^\infty(\mcs)$ is an FIO  in $I^{1/4}(
\mcn, \mcs; C_0)$ with canonical relation  
\beqq
C_0 = \{(z, \zeta, \bar z, \bar \zeta) \in T^* \mcn \backslash 0 \times T^*\mcs \backslash 0 : \bar z = z, \bar \zeta = \zeta|_{T_{\bar z}\mcs}\}
\eeqq
We consider the canonical relation $C_{wv}$ defined by
\beqq\label{eq-canowv}
\begin{gathered}
C_{wv} = \{ (w, \iota, \bar z, \bar \zeta) \in T^*\mcn \backslash 0 \times T^*\mcs \backslash 0: \text{$(w, \iota)$ is on the bicharacteristic }\\
\text{strip through some $(\bar z, \zeta)$ such that } \bar \zeta = \zeta|_{T_{\bar z}\mcs} \text{ and } \mcp(\bar z, \zeta) = 0\}
\end{gathered}
\eeqq
The  next result is straight forward from Theorem 5.1.2 of \cite{Dui}.
\begin{prop}\label{prop-wvpara}
There exists $E_1 \in I^{-1/4}(\mcn, \mcs; C_{wv}), E_2 \in  I^{-5/4}(\mcn, \mcs; C_{wv})$ such that 
\beqq\label{eq-parawv}
\begin{gathered}
P(z, D)E_k \in C^\infty(\mcn), \quad k = 1, 2\\
\rho_0 E_1 - \id \in C^\infty(\mcs), \quad \rho_0 E_2 \in C^\infty(\mcs)\\
\rho_0 D_t E_1 \in C^\infty(\mcs), \quad  \rho_0 D_t E_2  - \id \in C^\infty(\mcs)
\end{gathered}
\eeqq
\end{prop}
Now we can represent the solution of \eqref{eq-cauchy1} as $
u = E_1 f_1 + E_2 f_2
$ modulo a smooth term. Note that this is not the same representation used in Section \ref{sec-mod} because the Cauchy data are not re-parametrized to $h_1, h_2.$ 
It is natural to decompose $C_{wv}$ as the disjoint union of $C_{wv}^+$ and  $C_{wv}^-$ which are 
\beqq\label{eq-canowvpm}
\begin{gathered}
C^\pm_{wv} =  \{(w, \iota, \bar z, \bar \zeta) \in T^*\mcn \backslash 0 \times T^*\mcs \backslash 0: \text{$\iota$ is future/past }\\
\text{pointing light-like and lies on the bicharacteristic strip  through}\\
\text{ some $(\bar z, \zeta)$ such that }\bar \zeta = \zeta|_{T_{\bar z}\mcs} \text{ and } \mcp(\bar z, \zeta) = 0\}
\end{gathered}
\eeqq
These are  \eqref{eq-cano0} under the parametrization in Section \ref{sec-mod}. 
We can decompose (for $k = 1, 2$)
\beq
E_k = E_k^+ + E_k^-, \quad E_k^\pm \in I^{1- k - 1/4}(\mcn, \mcs; C_{wv}^\pm).
\eeq
We will find the relation of the principal symbols of $E_1^\pm, E_2^\pm$. We remark that the Maslov bundle and the half density bundle can be trivialized because the Lagrangians involved allow global parametrization. We will not show these factors in the notations below. 
\begin{lemma}\label{lm-wvsym}
Let $e_k^\pm, k = 1, 2$ be the principal symbol of $E_k^\pm$ on $\La^\pm = (C^{\pm}_{wv})'$ respectively. Suppose that the sub-principal symbol of $P(z, D)$ is purely imaginary, in which case $P(z, D)$  is of the form
\beqq\label{eq-hyper1}
P(z, D) = \square + \sum_{j = 0}^n \imath A_j(z)D_j + B(z)
\eeqq
where $A_j(z)$ are real valued smooth functions. Then $e_k^\pm, k = 1, 2$ are  real valued and 
\beq
e_1^+ > 0, \quad e_2^+ >0, \quad e_1^- >0, \quad e_2^- <0.
\eeq
\end{lemma}
\bpf
We  can find the principal symbols  following the argument in \cite[page 117]{Dui}. For $E_k \in I^{1-k-1/4}(\mcn, \mcs; C_{wv}), k = 1, 2$, if $e_k$ is the principal symbol of $E_k$, then it satisfies
\beqq\label{eq-lie}
\frac{1}{\imath} \mcl_{H_{p}}e_k + p_{sub} e_k = 0
\eeqq
where $\mcl_{H_p}$ denotes the Lie derivative and $p_{sub}$ denotes the sub-principal symbol of $P(z, D)$. This is a transport equation along the bicharacteristics. The initial conditions are determined as follows. For $(\bar z, \bar \zeta)\in T^*\mcs$, we have two cotangent vectors $(\bar z, \zeta^\pm)$ corresponding to it in $T^*\mcn$ where (regarding $\bar \zeta$ as a covector on $\mcs$)
\beq
\zeta^+ = (\tau, \bar\zeta), \quad \zeta^- = (-\tau, \bar \zeta)
\eeq
where $\tau = |\bar \zeta|$.  
From  the initial conditions in \eqref{eq-parawv}, we have
\beqq\label{eq-bicharini}
\begin{gathered}
e_1(\bar z, \zeta^+; \bar z, \bar \zeta)  = e_1(\bar z, \zeta^-; \bar z, \bar \zeta) >0  \\
e_2(\bar z, \zeta^+; \bar z, \bar \zeta) = - e_2(\bar z, \zeta^-; \bar z, \bar \zeta)   >0
\end{gathered}
\eeqq
which are all real valued. 
Let $\gamma^\pm(s), s \in \mbr$ be bicharactersitics such that $\gamma^\pm(0) = (\bar z, \zeta^\pm) \in T^*\mcn$. Along $\gamma^\pm(s)$, the equation \eqref{eq-lie} can be written as 
\beqq\label{eq-trans}
\p_s e_k(s) + a(s)e_k(s) = 0, \quad e_k(0) = e_k^\pm(\gamma^\pm(0))
\eeqq
where $a(s) = \imath p_{sub}(\gamma^\pm(s))$.  
Solving \eqref{eq-trans}, we obtain that 
\beqq\label{eq-bicharsol}
e_k(s) = e_k(0) e^{-\int_0^s a(\beta) d\beta}
\eeqq
Consider the operator \eqref{eq-hyper1}. 
In local coordinates, let $p_2(z, \zeta)$ be the symbol modulo $S^{0}(T^*\mcn)$, namely
\beq
p_2(z, \zeta) = g(\zeta, \zeta) +  \sum_{j = 0}^n \imath A_j(z)\zeta_j
\eeq
where $g(\zeta, \zeta) = -\tau^2 +|\xi|^2, \zeta = (\tau, \xi)$.   We have modulo symbols of order $0$
\beq
p_{sub}(z, \zeta) = \sum_{j = 0}^n \imath A_j(z)\zeta_j - \frac{1}{2\imath} \sum_{j = 0}^n \frac{\p^2 p_2(z, \zeta)}{\p z_j \p \zeta_j} = \sum_{j = 0}^n \imath A_j(z)\zeta_j   - \frac{1}{2\imath}\sum_{i, j = 1}^n\frac{\p^2 (h_{ij}\zeta^i\zeta^j)}{\p z_j \p \zeta_j}  
\eeq
If the subprincipal symbol $p_{sub}(z, \xi)$ is pure imaginary,  the coefficients of the transport equation \eqref{eq-trans} are real valued. We can tell from \eqref{eq-bicharsol} that $e_k^\pm$ are real valued and the signs are determined by the initial conditions in \eqref{eq-bicharini}. This completes the proof.  
\epf

\section{The composition as Lagrangian distributions}\label{sec-comp2}
In this section, we re-examine the composition in Section \ref{sec-comp1} from the point of view of Lagraigian distributions and complete the proof of Theorem \ref{thm-main}. Let's outline the main ingredients. 
We look at $E^*\chi N \chi_{[t_0, t_1]} E, \text{ with } E = E_k^\pm, k = 1, 2$ in Proposition \ref{prop-wvpara}.  
\begin{enumerate}[(i)]
\item As $\chi \cdot \chi_{[t_0, t_1]} = 0$, from Section \ref{sec-rep}, we know that $\chi N \chi_{[t_0, t_1]} \in I^{-n/2}(\mbr^{n+1}, \mbr^{n+1}; \La_1)$ at least when the characteristic function $\chi_{[t_0, t_1]}$ were smooth. Note that the role of $\chi$ is to keep the kernel of $N$ away from the diagonal $\La_0$ where the principal symbol is singular!
\item We will show that $\La_1$ intersect $\La_\pm$ cleanly with excess one so the composition $\chi N \chi_{[t_0, t_1]} E \in I^{\ast}(\mcn, \mcs; C_{wv})$ as  a result of Duistermaat-Guillemin's clean FIO calculus with the order $\ast$ to be determined. For this, we need to address some issue caused by the characteristic function. 
\item We can compose the operator in (ii) with $E^*$ by using clean FIO calculus again to conclude that $E^*\chi N \chi_{[t_0, t_1]} E \in \Psi^{\ast}(\mcs)$.  
\end{enumerate}
This is what behind the calculations in Section \ref{sec-comp1}. In the follows, we carry out the details of the above arguments. In this section, we assume that $n\geq 2$ is an integer.

\begin{lemma}\label{lm-clean}
Consider $\La_1$ defined in \eqref{eq-lag2}. Then $\La_1$ intersects $\La^\pm = (C_{wv}^\pm)'$ cleanly with excess one. 
\end{lemma}
\bpf
We check by the definition of clean intersection. We use the following parametrization for $\La_1$
 \beq 
 \begin{gathered}
 \La_1 = \{(t, x, \tau, \xi; t', x', \tau', \xi')\in T^*\mbr^{n+1}\backslash 0\times T^*\mbr^{n+1}\backslash 0: x = x' + (t - t')\xi/|\xi|,\\
 \tau = \pm|\xi|,  \tau' = -\tau,  \xi' = -\xi\}.
  \end{gathered}
 \eeq  
For $\La^\pm$, we use
 \beq
 \begin{gathered}
 \La^\pm = \{(\tilde t, \tilde x, \tilde \tau, \tilde \xi; z,  \eta) \in T^*\mbr^{n+1}\backslash 0 \times T^*\mbr^{n}\backslash 0 :  \tilde x = z \pm \tilde t\eta/|\eta|, \tilde \xi =- \eta, \tilde \tau = \pm |\eta|\}.
  \end{gathered}
 \eeq
We consider $\La^+$ below. The case for $\La^-$ is similar. Let $\mcx = \La_1\times \La^+$ and $\mcy = T^*\mcm\times \diag(T^*\mcm)\times T^*\mcs$. These are submanifolds of $T^*\mcm\times T^*\mcm\times T^*\mcm\times T^*\mcs.$ We show that for $p\in \mcx\cap \mcy$,  $T_p\mcx \cap T_p\mcy = T_p(\mcx\cap \mcy)$. 

First of all, $\La_1$ is parametrized by $ (t, t', x', \xi') \in \mca \doteq \mbr\times \mbr\times \mbr^n\times \mbr^n$. 
Also, $\La^+$ is parametrized by $(\tilde t, z, \eta) \in \mcb \doteq \mbr\times \mbr^n \times \mbr^n$ as 
\beq
\tilde x = z + \tilde t\eta/|\eta|, \quad \tilde \tau  = |\eta|, \quad \tilde \xi = -\eta.
\eeq
Consider $p  \in \mcx\cap \mcy$ and if we write $p = (t, t', x', \xi', \tilde t,   z,  \eta)\in \mca\times \mcb$, we must have 
\beqq\label{eq-q}
t' = \tilde t, \quad x' = z + \tilde t\eta/|\eta|, \quad |\eta| = |\xi'|, \quad \xi' = -\eta.
\eeqq
Thus, $\mcx\cap \mcy$ is parametrized by $ (t, \tilde t, z, \eta) \in \mcd \doteq \mbr\times \mbr\times \mbr^n\times \mbr^n$ as 
\beq
\begin{gathered}
t, \quad x = z +   t  \eta/|\eta|, \quad \tau  = -|\eta|,\quad \xi = \eta,   \quad t' = \tilde t, \quad x' = z + \tilde t \eta/|\eta|, \\
 \tau' =  |\eta|, \quad \xi' =  -\eta, \quad   \quad \tilde t, \quad  z, \quad \eta
\end{gathered}
\eeq
We find that tangent vector $\delta p \in T_p(\mcx\cap \mcy)$ is given by 
\beqq\label{eq-tempp1}
\begin{gathered}
\delta p =  (\delta t, \delta z + \eta/|\eta|\delta  t   +  t \alpha d\eta, -\eta/|\eta|\delta \eta, \delta \eta, \\
\delta \tilde t,  \delta z + \eta/|\eta|\delta \tilde t + \tilde t \alpha d\eta,  \eta/|\eta|\delta \eta, -\delta \eta,\\
\delta \tilde t, \delta z + \eta/|\eta|\delta \tilde t + \tilde t \alpha d\eta, \eta/|\eta|\delta \eta, -\delta \eta, \delta z, \delta \eta)
\end{gathered}
\eeqq
where  $\alpha = \p_\eta (\eta/|\eta|)$. 
Next, we compute $T_p \mcx$ and $T_p \mcy$ and find their intersection. For $\delta p \in T_p \mcx$, we use variables in $\mca$ and $\mcb$ to  get   
\beq
\begin{gathered}
\delta p  = (\delta t, \delta x' - (\delta t - \delta t')\xi'/|\xi'| - (t - t')\beta d\xi',  \mp \xi'/|\xi'|\delta \xi',-\delta \xi',  \delta t', \delta x', \pm \xi'/|\xi'|\delta \xi', \delta \xi', \\
\delta \tilde t, \delta z + \eta/|\eta|\delta \tilde t + \tilde t \alpha d\eta, \eta/|\eta|\delta \eta, -\delta \eta, \delta z, \delta \eta)
\end{gathered}
\eeq
where $\beta = \p_\eta(\eta/|\eta|)|_{\eta = \xi'}$. 
For $\delta p \in T_p \mcy$, we see that 
\beq
\delta t = \delta t', \quad \delta x' = \delta z + \eta/|\eta| \delta \tilde t + \tilde t \alpha d\eta, \quad \pm \xi'/|\xi'|\delta \xi' = \eta/|\eta|\delta \eta, \quad  \delta \eta = -\delta \xi'
\eeq
Also, at the intersection we use \eqref{eq-q}  and $(t, \tilde t, z, \eta)$ as variables to get 
\beqq\label{eq-tempp2}
\begin{gathered}
\delta p  = (\delta t, \delta z +   \eta/|\eta| \delta t + t \alpha d\eta, - \eta/|\eta|\delta \eta, \delta \eta,  \delta \tilde t, \delta z + \eta/|\eta| \delta \tilde t + \tilde t \alpha d\eta ,   \eta/|\eta|\delta \eta, -\delta \eta, \\
\delta \tilde t, \delta z + \eta/|\eta|\delta \tilde t + \tilde t \alpha d\eta, \eta/|\eta|\delta \eta, -\delta \eta, \delta z, \delta \eta)
\end{gathered}
\eeqq
Comparing \eqref{eq-tempp1} and \eqref{eq-tempp2}, we proved $T_p\mcx\cap T_p\mcy = T_p(\mcx\cap \mcy)$.

To find the excess, we see that $\codim(\mcx) = 8n+6  -(4n+3) = 4n + 3$, $\codim(\mcy) = (8n+6) -(6n+4) = 2n+2$. Also, $\dim(\mcx\cap \mcy) = 2n+2$. So the excess (see e.g. \cite[Appendix C.3]{Ho3})
\beq
e = \codim(\mcx) + \codim(\mcy) - \codim(\mcx\cap \mcy) = 4n+3 + 2n+2 - (6n+4) = 1.
\eeq 
This completes the proof of the lemma. 
\epf

\begin{lemma}\label{lm-comp3}
 The composition $\chi N \chi_{[t_0, t_1]} E_k^\pm \in I^{-n/2 + 1/4 + 1 - k}(\mcn, \mcs; C^\pm_{wv})$ and the principal symbol is non-vanishing. 
\end{lemma}
\bpf
First, we explain the idea by replacing  $\chi_{[t_0, t_1]}$ by a  smooth function $\tilde \chi$ compactly supported in $[t_0, t_1]$. Then because $\chi(t)\tilde \chi(t) = 0$, we know from Section \ref{sec-rep} that $\chi N \tilde \chi \in I^{-n/2}(\mcn, \mcn; \La_1)$. One can apply the clean calculus directly to see that $\chi N \tilde \chi E_k^\pm \in I^{-n/2 + 1/4 + 1 - k}(\mcn, \mcs; C^\pm_{wv})$ using Lemma \ref{lm-clean}.  For $p = (t, x, \tau, \xi, y, \eta)\in \La^\pm$, let $C_p$ be the fiber over $p$ in $T^*\mcm\times T^*\mcm \times T^*\mcs$ which is connected and compact. Then the principal symbol of the composition at $p$ is given by 
\beqq\label{eq-fiber}
\int_{C_p} \sigma(\chi N \tilde \chi)(t, x, \tau, \xi, t', x', \tau', \xi')\sigma(E_k^\pm)( t', x', \tau', \xi', y, \eta)
\eeqq
where $\sigma(\chi N \tilde \chi), \sigma(E_k^\pm)$ denote the principal symbols of $\chi N \tilde \chi, E_k^\pm$ respectively and the integration is over the fiber $C_p$, 
see \cite[Theorem 25.2.3]{Ho4}. Because  both symbols are real valued and non-vanishing on the fiber (modulo the Maslov factors), we see that the principal symbol of the composition is real valued and non-vanishing.

Next, for the characteristic function, the difference is that the fiber $C_p$ is connected but not compact. The arguments  in \cite[Theorem 25.2.3]{Ho4} still work provided that the integral \eqref{eq-fiber} is finite. To show this, we will extend the operators $E_k$ to a larger set. 

For $\eps > 0$ and small, we let $\tilde \mcn = (-\eps, T)\times \mbr^{n}$ which is an open set containing $\mcn$.   To extend $E_k$, we make use of the wave group property to re-parametrize the Cauchy data at   another Cauchy surface. We start with the backward Cauchy problem
\beq
\begin{gathered}
P(z, D)u(z) = 0, \quad z \in (-\eps, 0)\times \mbr^{n}\\
u = f_1,\quad D_t u = f_2, \quad \text{ at } \mcs.
\end{gathered}
\eeq
Let $B_1\in I^{-1/4}( (-\eps, 0)\times \mbr^n, \mcs; C_{wv}),$ and $B_2 \in I^{-5/4}( (-\eps, 0)\times \mbr^n, \mcs; C_{wv})$ be the parametrix of the backward Cauchy problem. Here, $C_{wv}^\pm$ is understood as the set \eqref{eq-canowvpm} but defined on $T^*((-\eps, 0)\times \mbr^n)\times T^*\mcs$.  Let $\rho_{\eps}$ be the restriction operator to $\mcs_{-\eps} = \{-\eps\}\times \mbr^{n}$. We let 
\beq
\begin{gathered}
f_1^{\eps} = \rho_{\eps} u = \rho_{\eps} B_1 f_1 + \rho_{\eps} B_2 f_2, \\
f_2^{\eps} = \rho_{\eps} D_t u = \rho_{\eps} D_t B_1 f_1 + \rho_{\eps} D_t B_2f_2
\end{gathered} 
\eeq
be the corresponding Cauchy data at $\mcs_{\eps}$. Now we consider the forward Cauchy problem 
\beq
\begin{gathered}
P(z, D)u(z) = 0, \quad z\in \tilde \mcn \\
u = f_1^{\eps},\quad D_t u = f_2^{\eps}, \quad \text{ at } \mcs_{-\eps}.
\end{gathered}
\eeq
Let $E_k^{\eps} \in I^{-1/4+1-k}(\tilde \mcn, \mcs_{-\eps}; C_{wv})$ be the corresponding parametrix by solving the forward Cauchy problem from $\mcs_{-\eps}.$ We consider 
\beq
\tilde u = E_1^{\eps} f_1^{\eps} + E_2^{\eps} f_2^{\eps}
\eeq
In particular, 
\beqq\label{eq-newu}
\begin{gathered}
\tilde u =  E_1^{\eps}(\rho_{\eps} B_1 f_1 + \rho_{\eps} B_2 f_2) + E_2^{\eps}( \rho_{\eps} D_t B_1 f_1 + \rho_{\eps} D_t B_2f_2)\\
 =  \tilde E_1 f_1 +    \tilde E_2 f_2
\end{gathered}
\eeqq
where 
\beq
\tilde E_1 = E_1^{\eps} \rho_{\eps} B_1 + E_2^{\eps} \rho_{\eps} D_t B_1, \quad \tilde E_2 = E_1^{\eps} \rho_{\eps} B_2 + E_2^{\eps} \rho_{\eps} D_t B_2
\eeq
First of all, $\rho_\eps B_k \in I^{1-k}(\mcs_{-\eps}, \mcs; C_\eps), k = 1, 2$ and $\rho_\eps D_t B_k\in  I^{2-k}(\mcs_{-\eps}, \mcs; C_\eps)$ are FIOs associated with canonical relations 
\beq
\begin{gathered}
C_\eps = \{(x, \xi, y, \eta)\in T^*\mcs_{-\eps}\backslash 0\times T^*\mcs \backslash 0: \text{$(x, \tilde \xi)$ and $(y, \tilde \eta)$ are on} \\
\text{the same bicharacteristics of $H_\mcp$ for } \tilde \xi|_{T_{x}\mcs_\eps} = \xi, \tilde \eta|_{T_y\mcs} = \eta\}
\end{gathered}
\eeq
see \cite[Chapter 5]{Dui}. 
We observe that the canonical relation is a canonical graph. For the composition $E^\eps_j \rho_s B_k$, we can use the transversal intersection calculus to conclude that 
\beq
E^\eps_j \rho_s B_k \in I^{-1/4 +1 - j + 1 - k}(\tilde \mcn, \mcs; C_{wv}) = I^{-1/4 +2 - j - k}( \tilde \mcn,  \mcs; C_{wv}).
\eeq
Thus, 
\beq
\tilde E_1 \in I^{-1/4}(\tilde \mcn, \mcs; C_{wv}), \quad \tilde E_2 \in I^{-5/4}(\tilde \mcn, \mcs; C_{wv})
\eeq
By the regularity of Cauchy problem,  these are extensions of $E_1, E_2$ which means that the  Schwartz kernel of $\tilde E_k - E_k$ belongs to $C^\infty(\mcn \times \mcs), k = 1, 2.$

Let $\psi(t), t\in \mbr$ be a smooth cut-off function such that $\psi(t) = 1$ if $t\in [-\eps/2, t_1+\eps/2]$ and $\psi(t) = 0$ if $t\notin (-\eps, t_1 +\eps)$. We see that the Schwartz kernels of $\psi \tilde E_k, k = 1, 2$ are equal to those of $E_k$ on $ \mcn \times \mcs$ and $ \psi \tilde E_k$ are compactly supported.  Then we consider the composition $\chi N  \psi \tilde E_k, k = 1, 2$. We can now apply the clean FIO calculus Theorem 25.2.3  of \cite{Ho4} to get that $\chi N  \psi \tilde E_k$ are Fourier integral operators in $I^{-n/2 + 1/4 + 1 -k}(\tilde N, \mcs, C_{wv})$. We can use \cite[Proposition 25.1.5']{Ho4} and conclude that the principal symbol of $\chi N \psi \tilde E_k$ is given by  
\beqq\label{eq-inte1} 
\int \sigma(\chi N\tilde \psi) \sigma(\psi \tilde E_k)
\eeqq
where $\tilde \psi$ is another compactly smooth cut-off function such that $\tilde \psi\psi = \psi$.  The integration is over the fiber $\gamma(t), t\in [-\eps, t_1 +\eps].$ Over $[0, t_1]$, the integral is well-defined. This shows that the composition $\chi N \chi_{[0, t_1]}E_k\in I^{-n/2 + 1/4 + 1 - k}(\tilde N, \mcs, C_{wv})$. 
\epf

\begin{lemma}\label{lm-comp4}
For $j, k = 1, 2$, we have
\begin{enumerate}
\item $E^{\pm, \ast}_j \chi N \chi_{[t_0, t_1]} E_k^\pm  \in \Psi^{-n/2+ 1/2 + 2 - j - k}(\mcs)$ are elliptic. 
\item $E^{+, \ast}_j \chi N \chi_{[t_0, t_1]} E^-_k, E^{-, \ast}_j \chi N \chi_{[t_0, t_1]} E^+_k$ are smoothing operators on $\mcs.$
\end{enumerate}
\end{lemma}
\bpf
 First of all, $E^{\pm, *}_j \in I^{-1/4 + 1-j}(\mcs, \mcm; C_{wv}^{\pm, -1})$ and $\chi N \chi_{[t_0, t_1]} E_k^\pm \in I^{-n/2 + 1/4 + 1-k}$ $(\mcm, \mcs; C_{wv}^\pm)$. Let $\La^\pm = (C_{wv}^\pm)'$ and $\La^{\pm, -1} = (C_{wv}^{\pm, -1})'$. We first prove that $\La^{\pm, -1}$ intersect $\La^\pm$ cleanly with excess one. 
 
We consider the plus sign. Recall that 
  \beqq
 \begin{gathered}
 \La^+ = \{( t,  x,   \tau,   \xi; z,  \eta) \in T^*\mbr^{n+1}\backslash 0 \times T^*\mbr^{n}\backslash 0 :  x = z +  t\eta/|\eta|,   \xi =- \eta,  \tau =  |\eta|\} 
  \end{gathered}
 \eeqq
 and it can be parametrized by $(t, z, \eta) \in \mcb = \mbr\times \mbr^n \times \mbr^n$. Let $\mcx = \La^{+, -1}\times \La^+$ and $\mcy = T^*\mcs\times \diag(T^*\mcm)\times T^*\mcs$. We see that 
 \beq
 \begin{gathered}
 \mcx\cap \mcy = \{(\tilde z, \tilde \eta, \tilde t,  \tilde x,   \tilde \tau,   \tilde \xi;  t,  x,   \tau,   \xi, z,  \eta): x = z +  t\eta/|\eta|,   \xi =- \eta,  \tau =  |\eta|\\
\eta = \tilde \eta, \tau = \tilde \tau,  \xi = \tilde \xi,  t = \tilde t, \tilde x = x, \tilde z = z\}
 \end{gathered}
 \eeq
 So this set is parametrized by $t, z, \eta$. For $p\in \mcx\cap \mcy$, the tangent vector $\delta p\in T_p(\mcx\cap \mcy)$ is 
 \beqq\label{eq-q1}
 \begin{gathered}
 \delta p = (\delta z, \delta \eta, \delta t, \delta z + \eta/|\eta| \delta t + t \p_\eta(\eta/|\eta|) \delta \eta, \eta/|\eta|\delta \eta, -\delta \eta\\
\delta t, \delta z + \eta/|\eta| \delta t + t \p_\eta(\eta/|\eta|) \delta \eta, \eta/|\eta|\delta \eta, -\delta \eta, \delta z, \delta \eta)
 \end{gathered}
 \eeqq
 Next, for $p \in \mcx$ which is parametrized by $(t, z, \eta, \tilde t, \tilde z, \tilde \eta)$, the tangent vector is given by 
  \beqq\label{eq-q2}
 \begin{gathered}
 \delta p = (\delta \tilde z, \delta \tilde\eta, \delta \tilde t, \delta \tilde z + \tilde\eta/|\tilde\eta| \delta \tilde t + \tilde t \p_{\tilde\eta}(\tilde\eta/|\tilde\eta|) \delta \tilde\eta, \tilde\eta/|\tilde\eta|\delta \tilde\eta, -\delta \tilde \eta\\
\delta t, \delta z + \eta/|\eta| \delta t + t \p_\eta(\eta/|\eta|) \delta \eta, \eta/|\eta|\delta \eta, -\delta \eta, \delta z, \delta \eta)
 \end{gathered}
 \eeqq
 If $\delta p$ also belongs to $T_p \mcy$, we see that 
 \beq
 \delta t = \delta \tilde t, \quad \delta x = \delta\tilde x, \quad \delta \tau = \delta \tilde \tau, \quad \delta \eta = \delta \tilde \eta
 \eeq
 which implies $\delta z = \delta \tilde z$,  so \eqref{eq-q2} agree with \eqref{eq-q1}. This shows the intersection is clean. To find the excess, we see that $\text{codim}(T_p \mcx) =(8n + 4) - (4n + 2) = 4n + 2$ and $\text{codim}(T_p\mcy) =(8n + 4) -( 6n + 2) = 2n + 2$. Also $\text{codim}(T_p(\mcx\cap \mcy)) =(8n + 4) -(2n + 1) = 6n + 3$. Thus the excess 
 \beq
 e = 4n+2 + 2n + 2- 6n + 3 = 1. 
 \eeq
 
Now we can use the clean FIO calculus \cite[Theorem 25.2.3]{Ho4} to conclude that $E^{\pm, \ast}_j \chi N \chi_{[t_0, t_1]} E_k^\pm$ $ \in \Psi^{-n/2+ 1/2 + 2 - j - k}(\mcs)$. As both principal symbols of $E^{\pm, \ast}_j$ and $\chi N \chi_{[t_0, t_1]} E_k^\pm$ are real and non-vanishing (modulo the Maslov factors), the principal of the composition is the integration of the product of principal symbols so is also non-vanishing. This proves part (1). 

Part (2) can be seen from a wave front set analysis using e.g.\ \cite[Theorem 1.3.7]{Dui}. 
 \epf

 \bpf[Proof of Theorem \ref{thm-main}]
The idea is similar to that for Theorem \ref{thm-mainsim}, despite that the parametrization of the Cauchy data is different. We write the solution of \eqref{eq-cauchy1} as 
\beq
u = E_1^+ f_1 + E_2^+ f_2  + E_1^- f_1 +  E_2^-f_2.
\eeq 
Next, we apply $\chi L^*$ to $Lu$ to get  
\beqq\label{eq-u2}
\begin{gathered}
\chi N u = \chi N \chi_{[t_0, t_1]} E_1^+ f_1 + \chi N \chi_{[t_0, t_1]} E_2^+ f_2 
 + \chi N \chi_{[t_0, t_1]} E_1^- f_1 +  \chi N \chi_{[t_0, t_1]} E_2^-f_2.
 \end{gathered}
\eeqq
Now we apply $E^{+, \ast}_{1}$ and use part (2) of Lemma \ref{lm-comp4} to get 
\beqq\label{eq-upara1}
\begin{gathered}
E_1^{+, *}\chi N u = E_1^{+, *} \chi N \chi_{[t_0, t_1]} E_1^+ f_1 + E_1^{+, *} \chi N \chi_{[t_0, t_1]} E_2^+ f_2 + R_1 f_1 +  R_2 f_2
 \end{gathered}
\eeqq
with $R_1, R_2$ smoothing operators. In the following, we use $R_1, R_2$ to denote generic smoothing operators which may change line by line. From Lemma \ref{lm-comp4} part (1), we see that $E_1^{+, *} \chi N \chi_{[t_0, t_1]} E_1^+ \in \Psi^{-n/2 + 1/2}(\mcs)$ and $E_1^{+, *} \chi N \chi_{[t_0, t_1]} E_2^+ \in \Psi^{-n/2 -1/2}(\mcs)$. 

On the other hand, we apply $E_1^{-, \ast}$ to \eqref{eq-u2} to get 
\beqq\label{eq-upara2}
\begin{gathered}
E_1^{-, *}\chi N u = E_1^{-, *} \chi N \chi_{[t_0, t_1]} E_1^- f_1 + E_1^{-, *} \chi N \chi_{[t_0, t_1]} E_2^- f_2 + R_1 f_1 +  R_2 f_2
 \end{gathered}
\eeqq
From Lemma \ref{lm-comp4} part (1), we see that $E_1^{-, *} \chi N \chi_{[t_0, t_1]} E_1^- \in \Psi^{-n/2 + 1/2}(\mcs)$ and $E_1^{-, *} \chi N \chi_{[t_0, t_1]} E_2^- \in \Psi^{-n/2 -1/2}(\mcs)$. It follows from Lemma \ref{lm-wvsym} and the composition results Lemma \ref{lm-comp3} and \ref{lm-comp4} that 
\beq
\begin{gathered}
\sigma(E_1^{+, *} \chi N \chi_{[t_0, t_1]} E_1^+) > 0, \quad \sigma( E_1^{-, *} \chi N \chi_{[t_0, t_1]} E_1^-) > 0\\
\sigma(E_1^{+, *} \chi N \chi_{[t_0, t_1]} E_2^+ ) >0,  \quad \sigma( E_1^{-, *} \chi N \chi_{[t_0, t_1]} E_2^-) <0. 
\end{gathered}
\eeq
Let $Q^+, Q^- \in \Psi^{n/2 - 1/2}(\mcs)$ be parametrices for $E_1^{+, *} \chi N \chi_{[t_0, t_1]} E_1^+, E_1^{-, *} \chi N \chi_{[t_0, t_1]} E_1^-$ respectively. We know that the principal symbols of $Q^\pm$ are positive. Applying $Q^\pm$ to \eqref{eq-upara1}, \eqref{eq-upara2}, we get 
\beqq\label{eq-upara3}
\begin{gathered}
Q^+ E_1^{+, *}\chi N u =   f_1 + B_+ f_2 + R_1 f_1 +  R_2 f_2 
\end{gathered}
\eeqq
\beqq\label{eq-upara4}
\begin{gathered} 
Q^- E_1^{-, *}\chi N u =  f_1 + B_- f_2 + R_1 f_1 +  R_2 f_2
\end{gathered}
\eeqq
where
\beq
B_+ = Q_+E_1^{+, *} \chi N \chi_{[t_0, t_1]} E_2^+, \quad B_- = Q_-E_1^{-, *} \chi N \chi_{[t_0, t_1]} E_2^- 
\eeq
From \eqref{eq-upara3}, we get 
\beq
Q^+ E_1^{+, *}\chi N u  - Q^- E_1^{-, *}\chi N u =  (B_+ - B_-) f_2 + R_1  f_1 +  R_2  f_2
\eeq
Note that $B_\pm \in \Psi^{-1}(\mcs)$ are elliptic. Also, the principal symbol of $B_+$ is positive but the principal symbol of $B_-$ is negative. Thus $B_+ - B_- \in \Psi^{-1}(\mcs)$ is elliptic. Let $W\in \Psi^1(\mcs)$ be a parametrix for $B_+ - B_-$. We get 
\beqq\label{eq-fpara1}
WQ^+ E_1^{+, *}\chi N u  - WQ^- E_1^{-, *}\chi N u =  f_2 + R_1  f_1 +  R_2  f_2
\eeqq
So we solved $f_2$ up to smooth terms. We can use $f_2$ for example in \eqref{eq-upara3} to get  
\beqq\label{eq-fpara2}
Q^+ E_1^{+, *}\chi N u - B_+(WQ^+ E_1^{+, *}\chi N u  - WQ^- E_1^{-, *}\chi N u) =   f_1  + R_1 f_1 +  R_2 f_2\\
\eeqq

From this point, we can follow the proof of Theorem \ref{thm-mainsim} line by line. In fact, $WQ^\pm \in \Psi^{n/2 + 1/2}(\mcs)$ and $Q^+, B_+WQ^\pm \in \Psi^{n/2- 1/2}(\mcs)$ so 
\beq
\begin{gathered}
WQ^\pm: H_{\comp}^s(\mcs) \rightarrow H_{\loc}^{s - n/2 - 1/2}(\mcs)\\
Q^+, B_+WQ^\pm:  H_{\comp}^s(\mcs) \rightarrow H_{\loc}^{s - n/2 + 1/2}(\mcs)
\end{gathered}
\eeq
are bounded. 
For $E_k^{\pm, *}, k = 1, 2$, we have $E_k^{\pm, \ast}: H_{\comp}^{s}(\mcn) \rightarrow  H_{\loc}^{s + 1 - k}(\mcs)$ is bounded. Finally, $L^*: H_{\comp}^{s}(\mcc)\rightarrow H_{\loc}^{s+ 1/2}(\mbr^{n+1})$ is bounded for $n\geq 3$, see  Proposition \ref{prop-sobo}. We obtain that 
\beq
E_1^{\pm, *} \chi L^* :   H_{\comp}^{s }(\mcc)\rightarrow H_{\loc}^{s + 1/2}(\mcn)
\eeq
is bounded. 
Thus using \eqref{eq-fpara1}, \eqref{eq-fpara2}, we get 
\beq
\begin{gathered}
\|f_1\|_{H^{s+1}(\mcs) } \leq C\| L u\|_{H^{s + n/2}(\mcc) } + C_\rho \|f_1\|_{H^{s - \rho}(\mcs) } +   C_\rho \|f_2\|_{H^{s - \rho}(\mcs) }\\
\|f_2\|_{H^s(\mcs) } \leq C\| L u\|_{H^{s + n/2}(\mcc) } + C_\rho \|f_1\|_{H^{s - \rho}(\mcs) } +   C_\rho \|f_2\|_{H^{s - \rho}(\mcs) }
\end{gathered}
\eeq
The rest of the proof are as in Theorem \ref{thm-mainsim}. 
 For $n=2$, one just need to use that $L^*: H_{\comp}^{s }(\mbr^2\times \mbs^{1}) \rightarrow H^{s+ 1/4}_{\loc}(\mbr^{3})$ is bounded from Proposition \ref{prop-sobo}.
 \epf

We  conclude the proof of Theorem \ref{thm-main} with two remarks. 
\begin{remark}
In the proof of Theorem \ref{thm-main}, we actually constructed operators $A_1, A_2$ such that 
\beq
A_1 L u = f_1 + R_1 f_1 + R_2 f_2, \quad A_2 Lu =  f_2 + R'_1 f_1 + R'_2 f_2
\eeq
where $R_1, R_2, R_1', R_2'$ are smoothing operators. The operators $A_1, A_2$ can be used to determine wave front set of $f_1, f_2$ from $Lu$. 
\end{remark}

\begin{remark}\label{remark2}
There are other ways to fine tune the normal operator $E^*L^*LE$ in \eqref{eq-normal0} to prove Theorem \ref{thm-main}.  We outline one possible construction and leave the details to interested readers. Instead of using $\chi$ compactly supported in $(t_1, T)$, we let $\rho_{\tilde T}$ be the restriction operator to $\mcs_{\tilde T} = \{\tilde T\}\times \mbr^n$ for some $\tilde T\in (t_1, T)$. In particular, $\rho_{\tilde T} \in I^{1/4}(\mcs_T, \mcm; C_0)$ in which 
\beq
C_0 = \{(y, \eta, t, x, \tau, \xi)\in T^*\mcs_{\tilde T}\backslash 0 \times T^*\mcm^\circ \backslash 0: y = x, \eta = \xi\}
\eeq
see (5.1.2) of \cite{Dui}. Recall from Lemma \ref{lm-comp3} that $\chi N \chi_{[t_0, t_1]}E_\pm \in  I^{-n/2 + 1/4}(\mcn, \mcs; C^\pm_{wv})$. Note that the composition $C^\pm_{wv, 0} = C_0\circ C^\pm_{wv}$ is given by 
\beq
C^\pm_{wv, 0} = \{(y, \eta, x, \xi) \in T^*\mcs_{\tilde T}\backslash 0\times T^*\mcs\backslash 0: y = x \pm \tilde T \xi, \eta = \xi\}
\eeq
which is a canonical graph. One can show that the composition is clean  as in Lemma \ref{lm-comp3} and obtain that 
\beqq\label{eq-rho}
\rho_{\tilde T} N \chi_{[t_0, t_1]}E_\pm \in  I^{-n/2 + 1}(\mcs_{\tilde T}, \mcs; C^\pm_{wv, 0}).
\eeqq
In particular, the principal symbol is non-vanishing. Now, \eqref{eq-rho} is an elliptic FIO of canonical graph type. We can find parametrix $Q_\pm \in I^{n/2-1}(\mcs, \mcs_{\tilde T}; C^{\pm, -1}_{wv, 0})$ such that 
\beq
Q_\pm \rho_{\tilde T} N \chi_{[t_0, t_1]}E_\pm = \id + R_\pm
\eeq 
where $R_\pm$ are smoothing operators. The rest of the argument goes as in the proof of Theorem \ref{thm-main}. 
 \end{remark}
 

\section{The source problem}\label{sec-source}
In this section, we consider the source problem 
\beqq\label{eq-source}
\begin{gathered}
P(z, D) u = f, \text{ on } \mcm\\
u = 0 \text{ for } t < t_0
\end{gathered}
\eeqq
where $f$ is compactly supported in $\mcm$. 
Consider the determination of $f$ from $Lu$. Before stating the main result, we explain the difference to the Cauchy problem.

According to  \cite{MeUh}, there exists a parametrix $E$ for \eqref{eq-source} such that $P(z, D)E = \id$ modulo a smoothing operator. The Schwartz kernel of $E$ belongs to $I^{-3/2, -1/2}(\mcm\times \mcm; \La_0, \La_1)$. It suffices to look at $L E f$. It is natural to apply $L^*$ and study $L^*L Ef = N Ef$. The Schwartz kernel of $N$ belongs to $I^{-n/2, n/2 -1} (\mcm\times \mcm; \La_0, \La_1)$ so both $N$ and $E$ are paired Lagrangian distributions of the flow-out type. One can apply the composition result  in  \cite{AnUh} to conclude that $NE\in I^{-n/2-1, n/2 -2} (\mcm\times \mcm; \La_0, \La_1)$. It is possible to find a parametrix for $NE$ within the class of paired Lagrangian distributions, however  the remainder term belongs to $I^\mu(\mcm\times \mcm; \La_1)$ for some $\mu\in \mbr$ rather than smooth, see \cite{Wan, GrUh0}. Moreover, although the parametrix is good for reconstructing space-like singularities, time-like singularities are lost and it is not clear whether one can determine light-like singularities of $f$.  
Below, we will assume that $\WF(f)$ is contained in $\Gamma^{sp}$ and use the kernel of $NE$ on $\La_0\backslash \La_1$ to stably determine $f$.  We remark that in general relativity, space-like singularities correspond to particles moving slower than the speed of light.

For $\delta>0$,  let $\Gamma^{sp}_\delta = \{(t, x, \tau, \xi)\in T^*\mcm: \tau^2 - |\xi|^2 >\delta\}$. 
 \begin{theorem}\label{thm-main1}
 Suppose that $f \in \mce'(\mcm)$ is supported in a compact set $\mcv$ of $\mcm$ and that $\WF(f) \subset \Gamma^{sp}_\delta$ for some $\delta>0$. Let $u$ be the solution of \eqref{eq-source}. Then there exists $C >0$ such that 
\beq
\| f\|_{H^{s}(\mcm)}\leq C \|L u\|_{H^{s + 3 -\frac{s_0}2}(\mcc)}
\eeq
with $s_0$ in Proposition \ref{prop-sobo} and $s\geq 0$.
\end{theorem}
 \bpf 
Because $\WF(f) \subset \Gamma^{sp}_\delta$, there exists an elliptic pseudo-differential operator $\chi(D)\in \Psi^0(\mcm)$ whose  symbol $\chi(x, \xi)$ is supported in $\Gamma^{sp}_{\delta/2}$ such that $\chi(D)f = f$ modulo a smooth term.  Thus, $N E f = N E \chi(D)f$ modulo a smooth term. Because $E\in I^{-3/2, -1/2}(\mcm\times \mcm; \La_0, \La_1)$, we claim that $E\chi(D)\in \Psi^{-2}(\mcm)$ 
with principal symbol $\chi(\tau, \xi)\sigma_0(E)(\tau, \xi)$ supported in $\Gamma^{sp}_{\delta/2}$. Here, $\sigma_0(E)$ denotes the principal symbol of $E$ on $\La_0.$ To see this, we can split $E = E_0 + E_1$ such that  $E_0\in \Psi^{-2}(\mcm)$ and $\WF(E_1)$ is sufficiently close to $\La_1.$ Then we know that $E_0 \chi(D)\in \Psi^{-3}(\mcm)$ and $E_1\chi(D)$ is a smoothing operator as the result of a wave front analysis using e.g.\ Theorem 1.3.7 of \cite{Dui} because the symbol of $\chi(D)$ is supported away from $\La_1.$  

It follows from the same argument that $NE\chi(D)\in \Psi^{-3}(\mcm)$ with principal symbol $\chi(\tau, \xi)\sigma_0(E)(\tau, \xi)\sigma_0(N)(\tau, \xi)$ which is non-vanishing. Thus, we can find a parametrix $Q \in \Psi^{3}(\mcm)$ of $NE\chi(D)$ such that 
\beq
Q NE\chi(D)= \id + R
\eeq
where $R$ is a smoothing operator. For $f\in \mce'(\mcm)$ with $\WF(f)\subset \Gamma^{sp}_\delta$, we actually have   $
Q NE = \id + R$ 
where we changed $R$ to another smoothing operator. Finally, we get that for any $\rho\in \mbr$, 
\beq
\|f\|_{H^s(\mcm)} \leq C\|N Ef\|_{H^{s + 3}(\mcm)} + C_\rho \|f\|_{H^{s-\rho}(\mcm)}
\eeq
for some $C, C_\rho >0.$ Using the estimate of $L^*$, we arrive at 
\beq 
\|f\|_{H^s(\mcm)} \leq C\|L E  f\|_{H^{s + 3 -\frac{s_0}2}(\mcc)} + C_\rho \|f\|_{H^{s-\rho}(\mcm)}
\eeq
with $s_0$ in Proposition \ref{prop-sobo}. Let $u$ be the solution of \eqref{eq-source} with source $f$. We get 
\beqq\label{eq-estcomp}
\|f\|_{H^s(\mcm)} \leq C\|L u\|_{H^{s + 3 -\frac{s_0}2}(\mcc)} + C_\rho \|f\|_{H^{s-\rho}(\mcm)}
\eeqq

By using the injectivity of $L$ as in the proof of Theorem \ref{thm-mainsim}, we can get rid of the last term as in Theorem 1.1 of \cite{VaWa}. We denote by $H^s_\mcv(\mcm)$ the function space consisting of $f \in H^s(\mcm)$ supported in $\mcv.$ Then the inclusion of $H^s_\mcv(\mcm)$ into $H^{s-\rho}_\mcv(\mcm), \rho>0$ is compact. We claim that 
\beqq\label{eq-finalest}
\|f\|_{H^s(\mcm)} \leq C\|L u\|_{H^{s + 3 -\frac{s_0}{2}}(\mcc)}  
\eeqq
for $f$ with $\WF(f)\subset \Gamma_\delta^{sp}.$ 
We argue by contradiction. Assume the above is not true. We can get a sequence $f^{(j)}, j = 1, 2, \cdots$ with unit norm in $H_\mcv^s(\mcm)$ and $\WF(f^{(j)})\subset \Gamma_\delta^{sp}$ such that $L u^{(j)}$ goes to $0$ in $H^{s + 3 - \frac{s_0}2}(\mcc)$ where $u^{(j)}$ is the solution of \eqref{eq-source} with source $f^{(j)}.$  By \eqref{eq-estcomp}, we conclude that $1=\|f^{(j)}\|_{H^s(\mcm)} \leq C'_\rho \|f^{(j)}\|_{H^{s-\rho}(\mcm)}$ for some constant $C'_\rho.$ This gives a weak limit $f$ in $H^s(\mcm)$ along a subsequence, which thus converges strongly in $H^{s-\rho}(\mcm)$. Therefore,  $\|f\|_{H^{s-\rho}(\mcm)}$ is bounded below by $1/C'_\rho$, thus non-zero. Now we use  the regularity estimate of the source problem $\|u\|_{H^{s+1}(\mcm)}\leq C \|f\|_{H^s(\mcm)}$ to conclude that $L u = 0$ with $u$ the solution of \eqref{eq-source} with source $f$. By the injectivity of $L$ we get $u = 0$ which gives $f = 0$ from the equation \eqref{eq-source}. We reached a contradiction which means \eqref{eq-finalest} holds. This finishes the proof. 
 \epf


\end{document}